\documentclass[12pt, draft]{amsart}

\newcommand{\Z}{\mathbb Z}
\newcommand{\T}{\mathbb T}
\newcommand{\R}{\mathbb R}
\newcommand{\N}{\mathbb N}
\newcommand{\G}{\mathcal G}
\newcommand{\M}{\mathcal M}
\renewcommand{\L}{\mathcal L}

\newcommand{\be}{\begin{equation}}
\newcommand{\ee}{\end{equation}}
\newcommand{\bea}{\begin{eqnarray}}
\newcommand{\eea}{\end{eqnarray}}
\newcommand{\bes}{\begin{eqnarray*}}
\newcommand{\ees}{\end{eqnarray*}}

\newtheorem{theorem}{Theorem}[section]
\newtheorem{proposition}[theorem]{Proposition}
\newtheorem{lemma}[theorem]{Lemma}

\allowdisplaybreaks[3]
\numberwithin{equation}{section}

\begin{document}

\title[An equivalence relation on wavelets in higher dimensions]
{An equivalence relation on wavelets in higher dimensions associated with matrix dilations}
\author{Biswaranjan Behera} 
\address{Statistics and Mathematics Unit, Indian Statistical Institute, 203, B.\ T.\ Road,
Calcutta 700108, India}
\email{br\_behera@yahoo.com}
\date{}
\keywords{wavelet, MSF-wavelet, wavelet set, dilation matrix, translation invariance}
\subjclass[2000]{42C40}

\begin{abstract}
We introduce an equivalence relation on the set of single wavelets of $L^2(\R^n)$ associated with an arbitrary dilation matrix. The corresponding equivalence classes are characterized in terms of the support of the Fourier transform of wavelets and it is shown that each of these classes is non-empty.
\end{abstract}

\maketitle
\section{Introduction}
Let $n\geq 1$ be an integer. An $n\times n$ matrix $A$ will be called a {\it dilation matrix} for $\R^n$ if $A$ preserves the lattice $\Z^n$, that is, $A\Z^n\subset\Z^n$, and all eigenvalues of $A$ have absolute value greater than 1. These conditions imply that $A$ has integer entries and $|\det A|$ is an integer greater than 1. Let $B=A^t$, the transpose of $A$, and $a=|\det A|=|\det B|$.

Since $A\Z^n$ is a normal subgroup of the abelian group $\Z^n$, we can form the cosets of $A\Z^n$ in $\Z^n$. It is a well-known fact that the number of distinct cosets of $A\Z^n$ in $\Z^n$ is equal to $a= |\det A|$ (\cite{GM}, \cite{Woj}). A subset of $\Z^n$ which consists of exactly one element from each of these $a$ cosets will be called a {\it set of digits} for $A$. It is easy to see that if $K$ is a set of digits for $A$, then so is $K-\mu$, where $\mu\in K$. Therefore, we can assume, without loss of generality, that $0\in K$.

A finite set of functions $\Psi= \{\psi^1,\psi^2,\dots,\psi^L\}\subset L^2(\R^n)$ is said to be a {\it multiwavelet} of $L^2(\R^n)$ associated with the dilation $A$ if the system 
\begin{equation*}
\label{E.mwav}
\{\psi^l_{jk}:=a^{j/2}\psi^l(A^j\cdot-k):j\in\Z, k\in\Z^n, 1\leq l\leq L\}
\end{equation*}
forms an orthonormal basis for $L^2(\R^n)$. If there is a single function $\psi\in L^2(\R^n)$ such that 
\begin{equation*}
\label{E.wav}
\{\psi_{jk}:=a^{j/2}\psi(A^j\cdot-k):j\in\Z, k\in\Z^n\}
\end{equation*}
is an orthonormal basis for $L^2(\R^n)$, then we say $\psi$ to be an {\it $A$-wavelet}.

A.\ Calogero~\cite{Cal} provided a characterization of multiwavelets associated with a general dilation matrix. This result extends to the higher dimensional situation a well-known theorem for dyadic wavelets in $L^2(\R)$, which was proved independently by 
G.\ Gripenberg and X.\ Wang (see Chapter~7 of~\cite{HW} and reference cited there). Dai, Larson, and Speegle~(\cite{DLS}) proved that, associated with every dilation matrix, there exist wavelets $\psi$ of $L^2(\R^n)$ of the form $|\hat\psi|=\chi_K$, where $K$ is a measurable subset of $\R^n$. Such wavelets are called MSF ({\it minimally supported frequency}) $A$-wavelets and we shall call the corresponding set $K$ an {\it $A$-wavelet set}. The theorem of Calogero for single wavelets of $L^2(\R^n)$ is the following.

\begin{theorem}
\label{T.wavelet}
A function $\psi\in L^2(\R^n)$ is an $A$-wavelet for $L^2(\R^n)$ if and only if the following are satisfied:
\begin{itemize}
\item[(i)] $\|\psi\|_2=1$,
\item[(ii)] $\sum\limits_{j\in \Z} |\hat\psi(B^j\xi)|^2=1 
     \quad\mbox{for a.e.}~\xi\in\R^n$, 
\item[(iii)] $\sum\limits_{j\geq 0} \hat\psi(B^j\xi)            
     \overline{\hat\psi(B^j(\xi+2q\pi))}= 0 
     \quad\mbox{for a.e.}~\xi\in\R^n, q\in \Z^n\setminus B\Z^n$.
\end{itemize}
\end{theorem}

We use the following definition of Fourier transform:
\[
\hat f(\xi)=\int_{\R^n} f(x)e^{-i\langle \xi, x\rangle} dx, \quad\xi\in\R^n.
\]

The $A$-wavelet sets (equivalently, the MSF $A$-wavelets) admit a simple characterization in terms of two geometric conditions.

\begin{theorem} 
\label{T.wavset}
A set $K\subset\R^n$ is an $A$-wavelet set if and only if the following two conditions hold:
\begin{enumerate}
	\item[(i)] $\{K+2k\pi:k\in\Z^n\}$ is a partition of $\R^n$, and 
	\item[(ii)] $\{B^j K:j\in\Z\}$ is a partition of $\R^n$.
\end{enumerate}
\end{theorem}

Given an $A$-wavelet $\psi$, we define a sequence of closed subspaces of $L^2(\R^n)$ as follows:
\[
V_j=\overline{{\rm span}}\{\psi_{l,k}:\mbox{$l<j$},\mbox{$k\in\Z^n$}\},\quad j\in\Z.
\]
It is easy to verify that these subspaces satisfy the following properties:
\begin{enumerate}
\item[(i)] $V_j\subset V_{j+1}$ for all $j\in\Z$, 
\item[(ii)] $f\in V_j$ if and only if $f(A\cdot)\in V_{j+1}$
 for all $j\in\Z$, 
\item[(iii)] $\cup_{j\in\Z} V_j$ is dense in $L^2(\R^n)$, 
\item[(iv)] $\cap_{j\in\Z} V_j = \{0\}$, and
\item[(v)] $V_0$ is invariant under translation by elements of the group $\Z^n$.
\end{enumerate}

In view of the last property, it is natural to ask the question of the existence of other (larger) groups of translations under which the space $V_0$ remains invariant. For each non-negative integer $r$, we shall prove the existence of $A$-wavelets for which the corresponding space $V_0$ is invariant under translation by elements of the form $A^{-r}k$ for all $k\in\Z^n$.

For $y\in\R^n$, let $T_y$ be the (unitary) translation operator on $L^2(\R^n)$ defined by 
$T_y f(x) = f(x-y)$. We consider the following groups of translation operators:
\[
\G_r =\{T_{A^{-r}k}:k\in\Z^n\},\ r\in\N\cup\{0\}, \quad{\rm and}\quad
\G_{\infty}=\{T_y:y\in\R^n\}.
\]
 
Let $\mathcal G$ be a set of bounded linear operators on $L^2(\R^n)$ and $V$ a closed subspace of $L^2(\R^n)$. We say that $V$ is $\mathcal G$-invariant if $Tf\in V$ for every $f\in V$ and $T\in{\mathcal G}$.
 
Let $\L_r$ denote the collection of all $A$-wavelets such that the corresponding space $V_0$ is $\G_r$-invariant. It is clear that $\L_0$ is the set of all $A$-wavelets, and
\[
\L_{0}\supset\L_{1}\supset\L_{2}\supset\cdots\supset
\L_r\supset\L_{r+1}\supset\cdots\supset\L_{\infty}.
\]

These inclusions naturally give rise to an equivalence relation on the collection of $A$-wavelets. The equivalence classes are given by 
$\M_r = \L_r\setminus \L_{r+1}$, $r\geq 0$, and $\M_\infty = \L_\infty$. 
Therefore, $\M_r$, $r\geq 0$, consists of those $A$-wavelets for which $V_0$ is  $\G_r$-invariant but not $\G_{r+1}$-invariant.

This equivalence relation was defined in~\cite{Web} for dyadic wavelets in $L^2(\R)$, where the author characterized the equivalence classes. He also proved that $\M_r,\ r=0,1,2,3$, are non-empty. Later, in~\cite{BM},~\cite{SW}, examples of wavelets of $L^2(\R)$ were constructed for each of these equivalence classes, by different methods. This equivalence relation was extended by the author to the wavelets of the Hardy space $H^2(\R)$ in~\cite{Beh}, where it was also shown that the corresponding equivalence classes are non-empty.

In section~\ref{S.two} we characterize the equivalence classes in terms of the support of the Fourier transform of wavelets 
and in section~\ref{S.three} we prove that each of them is non-empty.
\section{Characterization of the equivalence classes}
\label{S.two}
For $j\in\Z$ and $y\in\R^n$, define the dilation operators $D_j$ and translation operators $T_y$ on $L^2(\R^n)$ as follows:
\[
D_jf(x)= a^{j/2}f(A^jx)\quad{\rm and}\quad T_yf(x)=f(x-y).
\]
In this notation we have $\psi_{jk}=D_j T_k\psi$. Observe that for $j\in\Z$ and $y\in\R^n$, we have
\bes
T_yD_j f(x) & = & D_j f(x-y)\\
&=& a^{j/2}f(A^j(x-y))\\
&=& T_{A^jy}a^{j/2}f(A^jx)\\
&=& D_jT_{A^jy}f(x)
\ees
Hence, we have the following commutation relation:
\begin{equation}
\label{E.comm}
T_y D_j = D_jT_{A^j y}, \quad j\in\Z~{\rm and~}y\in\R^n. 
\end{equation}

We define another sequence, $\{W_j:j\in\Z\}$, of closed subspaces of $L^2(\R^n)$ by
$W_j=V_{j+1}\ominus V_j$, the orthogonal complement of $V_j$ in $V_{j+1}$. Clearly,
\begin{equation}
\label{E.decom}
L^2(\R^n)=V_0\oplus\Bigl(\bigoplus_{j\geq 0}W_j\Bigr).
\end{equation}
Observe that $W_0=\overline{\rm span}\{\psi(\cdot-k):k\in\Z^n\}$. Therefore, the wavelet $\psi$ belongs to the space $W_0$. For this reason, first we characterize this space.

If $f\in W_0$, then there exists a sequence $\{c_k:k\in\Z^n\}\in l^2(\Z^n)$ such that $f=\sum_{k\in\Z^n}c_k\psi(\cdot-k)$. Taking Fourier transform of both sides, we get $\hat f(\xi)=m(\xi)\hat\psi(\xi)$ for some $2\pi\Z^n$-periodic function $m\in L^2(\T^n)$, where $\T^n=[-\pi,\pi]^n$. Conversely, if $\hat f$ is of this form, then by considering the Fourier series of $m$ and reversing the above steps, we see that $f\in W_0$. Therefore, we have the following characterization of $W_0$.
\begin{equation}
\label{E.w0}
W_0=\{ f\in L^2(\R^n):\hat f(\xi)=m(\xi)\hat\psi(\xi)~\mbox{for some}~m\in L^2(\T^n)\}.
\end{equation}

We begin with a few lemmas.

\begin{lemma}
\label{L.vpvq}
Fix a non-negative integer $r$, and let $p$, $q$ be integers with $p<q$. Then the space $V_p$ (resp.\ $W_p$) is $\G_r$-invariant if and only if $V_q$ (resp.\ $W_q$) is $\G_{r+q-p}$-invariant.
\end{lemma}

\proof
By definition, $f\in V_p$ if and only if $D_{q-p}\in V_q$. Suppose that $g\in V_q$ and find $f\in V_p$ such that $D_{q-p}f=g$. Now, from the commutation relation~(\ref{E.comm}), we have for $k\in\Z^n$
\[
T_{A^{-(r+q-p)}k}D_{q-p}f(x)=D_{q-p}T_{A^{-r}k}f(x).
\]
The lemma now follows from this equation.
\qed

The structure of the space $V_0$ can be quite complicated; in some cases it may not even be generated by the $\Z^n$-translates of a finite number of functions. On the other hand, $W_0$ is generated by the $\Z^n$-translates of a single function, namely, the wavelet $\psi$. The following lemma allows us to work with $W_0$ instead of $V_0$.

\begin{lemma}
\label{L.v0w0}
The space $V_0$ is $\G_r$-invariant if and only if $W_0$ is $\G_r$-invariant.
\end{lemma}

\proof
Suppose that $V_0$ is $\G_r$-invariant. By Lemma~\ref{L.vpvq}, $V_1$ is   
$\G_{r+1}$-invariant. In particular, $V_1$ is $\G_r$-invariant. But $V_1=V_0\oplus W_0$. Hence, $W_0$ is $\G_r$-invariant.

Conversely, let $W_0$ is $\G_r$-invariant. Then, again by Lemma~\ref{L.vpvq}, $W_j$ is $\G_{r+j}$-invariant for all $j\geq 0$. In particular, $W_j$ is $\G_r$-invariant, hence so is $\oplus_{j\geq 0}W_j$. The $\G_r$-invariance of $V_0$ now follows from the decomposition (\ref{E.decom}).
\qed

\begin{lemma}
\label{L.matrix}
Let $M$ be a dilation matrix and $K_M$ be a set of digits for $M$. Then
\[
\sum\limits_{\mu\in K_M} e^{-2\pi i \langle M^{-1}\mu,\nu\rangle}=|\det M|\delta_{0\nu},
\quad {\it for}~\nu\in K_{M^t},
\]
where $K_{M^t}$ is a set of digits for the dilation $M^t$.
\end{lemma}

\proof
This is the orthogonal relation for the characters of the finite group $\Z^n/M\Z^n$ (see \cite{Rud}). Observe that the mapping
\[
\mu+M\Z^n \mapsto e^{-2\pi i\langle M^{-1}\mu,\nu\rangle},
 \quad\nu\in K_{M^t}
\]
is a character of the (finite) coset group $\Z^n/M\Z^n$. If
$\nu=0$, then there is nothing to prove.
Suppose that $\nu\not=0$, then there exists a $\mu'\in K_M$ such
that $e^{-2\pi i\langle M^{-1}\mu',\nu\rangle}\not=1$. Since $K_M$
is a set of digits for $M$, so is $K_M-\mu'$. Hence,
\begin{equation}
\label{E.digit}
 \sum\limits_{\mu\in K_M} e^{-2\pi i\langle M^{-1}(\mu-\mu'),\nu\rangle}=
 \sum\limits_{\mu\in K_M} e^{-2\pi i\langle M^{-1}\mu,\nu\rangle}.
\end{equation}
Now
\bes
 \sum\limits_{\mu\in K_M} e^{-2\pi i\langle M^{-1}\mu,\nu\rangle}
 & = & e^{-2\pi i\langle M^{-1}\mu',\nu\rangle}
       \cdot \sum\limits_{\mu\in K_M} e^{-2\pi i\langle M^{-1}(\mu-\mu'),\nu\rangle} \\
 & = &  e^{-2\pi i\langle M^{-1}\mu',\nu\rangle}
       \cdot \sum\limits_{\mu\in K_M} e^{-2\pi i\langle M^{-1}\mu, \nu\rangle},
       \quad{\rm by}~(\ref{E.digit}).
\ees
Therefore,
\[
\sum\limits_{\mu\in K_M} e^{-2\pi i\langle M^{-1}\mu,\nu\rangle }= 0,
\quad{\rm since}~e^{-2\pi i\langle M^{-1} \mu',\nu\rangle }\not=1.
\]
\qed

Applying Lemma~\ref{L.matrix} to the dilation matrix $A^r$, we get the following simple but useful fact.

\begin{lemma}
\label{L.lattice}
Fix $r\in\N$ and let $l\in\Z^n\setminus B^r\Z^n$. Then there exists $k\in\Z^n$ such that $\langle B^{-r}l,k\rangle \not\in\Z$.
\end{lemma}

\proof
Since $A$ is a dilation matrix, so is $A^r$. From Lemma~\ref{L.matrix}, we have 
\[
\sum_{\mu\in K_{A^r}} e^{-2\pi i\langle A^{-r}\mu,\nu\rangle}=0, 
\quad\nu\in K_{B^r}\setminus\{0\}.
\]
If $l\in\Z^n\setminus B^r\Z^n$, then $l=B^rp+\nu'$, where $p\in\Z^n$ and $\nu'\in K_{B^r}\setminus\{0\}$. Therefore, for $\mu\in K_{A^r}$, we have 
$\langle A^{-r}\mu,l\rangle  = \langle A^{-r}\mu,B^rp+\nu'\rangle  = \langle \mu,p\rangle +\langle A^{-r}\mu,\nu'\rangle $; hence,
\begin{equation}
\label{E.sum}
\sum\limits_{\mu\in K_{A^r}} e^{-2\pi i\langle A^{-r}\mu,l\rangle }= 
\sum\limits_{\mu\in K_{A^r}} e^{-2\pi i\langle A^{-r}\mu,\nu'\rangle }=0.
\end{equation}
Now, if $\langle B^{-r}l,k\rangle =\langle A^{-r}k,l\rangle \in\Z$ for all $k\in\Z^n$, then each term in the sum in (\ref{E.sum}) would be equal to 1, a contradiction.
\qed

Let $U$ and $V$ be two measurable subsets of $\R^n$. Then $U$ is said to be $2\pi\Z^n$-{\it translation equivalent} to $V$ if there exists a measurable partition $\{U_k:k\in\Z^n\}$ of $U$ such that $\{U_k+2k\pi:k\in\Z^n\}$ forms a partition of $V$.
A function $f$ on $\R^n$ is said to be $2\pi\Z^n$-{\it periodic} if $f(x+2m\pi)=f(x)$ for a.e.\ $x\in\R^n$ and all $m\in\Z^n$.

Let $\psi$ be an $A$-wavelet. Denote $E={\rm supp~}\hat\psi$. For $k\in\Z^n$, define
\[
E(k)=\{\xi\in E:\xi+2k\pi\in E\}=E\cap(E+2k\pi).
\]
\begin{theorem}
\label{T.lr}
Let $r\in\N$ and $\psi$ be an $A$-wavelet. Then $\psi\in\L_r$ if and only if $E(k)=\emptyset$ for every $k\in\Z^n\setminus B^r\Z^n$.
\end{theorem}

\proof
{\bf Necessity}. Let $\psi\in\L_r$. That is, $V_0$ is $\G_r$-invariant, hence by Lemma~\ref{L.v0w0}, $W_0$ is $\G_r$-invariant. Since $\psi\in W_0$, the $\G_r$-invariance of $W_0$ implies that $\psi(\cdot-A^{-r}k)\in W_0$ for all $k\in\Z^n$. Taking Fourier transform, we get $e^{-i\langle B^{-r}\cdot,k\rangle }\hat\psi\in\widehat{W}_0$ for all $k\in\Z^n$. By the characterization of $W_0$ (see~(\ref{E.w0})), this is equivalent to saying that for each $k\in\Z^n$, there exists a $2\pi\Z^n$-periodic function $\mu_k\in\ L^2(\T^n)$ such that $e^{-i\langle B^{-r}\xi,k\rangle }\hat\psi(\xi)=\mu_k(\xi)\hat\psi(\xi)$. Note that if $\xi\in E$ (so that $\hat\psi(\xi)\not=0$), then we have
\begin{equation}
\label{E.muk}
e^{-i\langle B^{-r}\xi,k\rangle }=\mu_k(\xi).
\end{equation}
Now suppose that there exists $l\in\Z^n\setminus B^r\Z^n$ such that $E(l)\not=\emptyset$. Hence, there exists a non-trivial set $F$ such that both $F$ and $F+2l\pi$ are subsets of $E$. Now, for a.e. $\xi\in F$ and for all $k\in\Z^n$, we have
\bes
& & e^{-i\langle B^{-r}(\xi+2l\pi),k\rangle }\hat\psi(\xi+2l\pi)
 = \mu_k(\xi+2l\pi)\hat\psi(\xi+2l\pi) \\
\Longrightarrow &  & e^{-i\langle B^{-r}\xi,k\rangle } 
e^{-2\pi i\langle B^{-r}l,k\rangle }\hat\psi(\xi+2l\pi)
= \mu_k(\xi)\hat\psi(\xi+2l\pi) \\
& & \hskip 3cm\quad({\rm since~}\mu_k~{\rm is}~2\pi\Z^n{\rm -periodic})\\
\Longrightarrow &  & e^{-2\pi i\langle B^{-r}l,k\rangle }
= 1, \quad({\rm by~(\ref{E.muk})}, {\rm and~since}~\xi+2l\pi\in E)\\
\Longrightarrow &  & \langle B^{-r}l,k\rangle \in\Z\quad {\rm for~all~} k\in\Z^n.
\ees
This is a contradiction to Lemma~\ref{L.lattice}.

{\bf Sufficiency}. Suppose that $E(k)=\emptyset$ for all $k\in\Z^n\setminus B^r\Z^n$. To establish $\psi\in\L_r$, it is sufficient to show that for all $k\in\Z^n$, we have 
\begin{equation}
\label{E.nukxi}
e^{-i\langle B^{-r}\xi,k\rangle }\hat\psi(\xi)=\nu_k(\xi)\hat\psi(\xi)
\end{equation}
for some $2\pi\Z^n$-periodic function $\nu_k\in L^2(\T^n)$. Define the map
\[
\tau:\R^n\longrightarrow\T^n,\quad \tau(x)=x+2m(x)\pi,
\]
where $m(x)$ is the unique element in $\Z^n$ such that $x+2m(x)\pi\in\T^n$. Since the relation $\sum_{m\in\Z^n}|\hat\psi(\xi+2m\pi)|^2=1$~a.e. holds for every $A$-wavelet, it is clear that $\tau:E\longrightarrow\T^n$ is onto. First, we want to find a subset $H$ of $E$ such that $\tau:H\longrightarrow\T^n$ is a bijection, that is, $H$ is $2\pi\Z^n$-translation equivalent to $\T^n$.

Fix $\xi\in\T^n$. Define $Z(\xi)=\{k\in\Z^n:\xi+2k\pi\in E\}$. Notice that for a.e. $\xi$, $Z(\xi)$ is non-empty. For $\xi\in E$, we choose $k_\xi=0$. Now, let $\xi\not\in E$. If there is an $l=(l_1,l_2,\dots,l_n)\in Z(\xi)$ such that $l_1>0$, then we choose $k_\xi=k$, where $k_1=\inf\{l_1:l=(l_1,l_2,\dots,l_n)\in Z(\xi)~{\rm and}~l_1>0\}$. Otherwise we choose $k_\xi=k$, where $k_1=\sup\{l_1:l=(l_1,l_2,\dots,l_n)\in Z(\xi)~{\rm and}~l_1<0\}$. Now define the set
\begin{equation}
\label{setf}
H=\{\xi+2k_\xi\pi:\xi\in\T^n\}.
\end{equation}

Observe that $H$ is $2\pi\Z^n$-translation equivalent to $\T^n$, by construction. For $k\in\Z^n$, we define the function $\nu_k$ on $H$ by 
\begin{equation}
\label{E.nuk}
\nu_k(\xi)=e^{-i\langle B^{-r}\xi,k\rangle}\quad(\xi\in H)
\end{equation}
and extend $2\pi\Z^n$-periodically to the whole of $\R^n$ so that (\ref{E.nukxi}) holds clearly on $H$. To complete the proof we need to show that (\ref{E.nukxi}) also holds on $E\setminus H$.

Now, for almost every $\xi\in E\setminus H$, there exists $\xi'\in H$ and $l_\xi\in\Z^n$ such that $\xi=\xi'+2l_\xi\pi$, since $H$ is $2\pi\Z^n$-translation equivalent to $\T^n$. That is, $E(l_\xi)\not=\emptyset$. Hence, by hypothesis, $l_\xi=B^rp_\xi$ for some $p_\xi\in \Z^n$. Therefore, we have
\bes
e^{-i\langle B^{-r}\xi,k\rangle }\hat\psi(\xi) & = & 
e^{-i\langle B^{-r}(\xi'+2B^rp_\xi\pi),k\rangle }
\hat\psi(\xi'+2l_\xi\pi) \\
& = & e^{-i\langle B^{-r}\xi',k\rangle }\hat\psi(\xi'+2l_\xi\pi) \\
& = & \nu_k(\xi')\hat\psi(\xi'+2l_\xi\pi),\quad({\rm by~(\ref{E.nuk}),~since}~\xi'\in H) \\
& = & \nu_k(\xi-2l_\xi\pi)\hat\psi(\xi) \\
& = & \nu_k(\xi)\hat\psi(\xi)\quad ({\rm since}~\nu_k~{\rm is}~2\pi\Z^n{\rm-periodic}). 
\ees
This completes the proof of the theorem.
\qed

Theorem~\ref{T.lr} allows us to characterize the equivalence classes in terms of the support of the Fourier transform of the wavelets.

\begin{theorem}
\label{T.char}
\begin{enumerate}
	\item[(a)]The equivalence class $\M_\infty$ is precisely the collection of all MSF $A$-wavelets.
	\item[(b)] An $A$-wavelet $\psi\in\M_r$, $r\geq 1$, if and only if $E(k)=\emptyset$ for all $k\in\Z^n\setminus B^r\Z^n$ but there exists $l\in\Z^n\setminus B^{r+1}\Z^n$ such that $E(l)\not=\emptyset$.
	\item[(c)] An $A$-wavelet $\psi\in\M_0$ if and only if there exists $l\in\Z^n\setminus B\Z^n$ such that $E(l)\not=\emptyset$.
\end{enumerate}
\end{theorem}

\proof
Item (a) is an easy generalization of the corresponding result for the one-dimensional dyadic wavelets, which is proved in~\cite{Web}. Item (b) follows from Theorem~\ref{T.lr}. It also follows from Theorem~\ref{T.lr} that an $A$-wavelet $\psi\in\L_1$ if and only if $E(l)=\emptyset$ for every $l\in\Z^n\setminus B\Z^n$. Since $\L_0$ consists of all $A$-wavelets, (c) follows from this and the fact that $\M_0=\L_0\setminus\L_1$.
\qed


\section{Construction of wavelets in $\M_r$}
\label{S.three}
Let $D$ be a subset of $\R^n$ satisfying (ii) in Theorem~\ref{T.wavset}. For $E\subset\R^n$, we define the {\it translation projection} $\tau$ and {\it dilation projection} $d$ as follows:
\[ \tau(E)=\bigcup_{k\in\Z^n}\{(E+2k\pi)\cap{\T^n}\}, {\rm and} \]
\[ d(E)=\bigcup_{j\in\Z}\{(B^jE)\cap D\}. \]

The following result gives sufficient conditions for a set $S\subset\R^n$ to be a subset of an $A$-wavelet set. We refer to~\cite{BS} for the proof.

\begin{theorem}
\label{T.subws}
Let $A$ be a dilation matrix, $B=A^t$, and  $S\subset\R^n$. Then $S$ is a subset of an $A$-wavelet set if the following conditions hold:
\begin{itemize}
	\item[(a)] $S\cap(S+2k\pi)=\emptyset$ for all $k\in\Z^n$,
	\item[(b)] $(B^jS)\cap S=\emptyset$ for all $j\in\Z$,
	\item[(c)] there exists $\epsilon>0$ such that $N_\epsilon(0)\cap\tau(S)=\emptyset$, where $N_\epsilon(x)$ is the ball of radius $\epsilon$ around $x$,
	\item[(d)] $D\setminus d(S)$ has a non-empty interior.
\end{itemize}
\end{theorem}

We shall also require another result proved in~\cite{BS} for expansive matrices, i.e., $n\times n$ matrices such that all eigenvalues have modulus greater than $1$. Hence, a dilation matrix is an expansive matrix which preserves the lattice $\Z^n$.

\begin{lemma}
\label{L.fixed}
Suppose $A$ is a dilation matrix, $k\in\Z^n\setminus\{0\}$, and $p\in\Z\setminus\{0\}$. For $\alpha\in\Z^n$ and $j\in\Z$, let $f_\alpha(x)=B^{-p}x-2\alpha\pi$ and $g_j(x)=B^{-j}x-2k\pi$. Let $Y$ be the collection of all fixed points of 
${\mathcal F}\cup{\mathcal G}$, where
\begin{equation}
\label{E.fg}
{\mathcal F}=\{f_\alpha:\alpha\in\Z^n\}\quad and\quad
{\mathcal G}=\{g_j:j\in\Z\}.
\end{equation}
Then for any $y\not\in \overline{Y}$, there exists $\epsilon>0$ such that 
\begin{equation}
\label{E.faga}
f_\alpha(N_\epsilon(y))\cap N_\epsilon(y)=g_j(N_\epsilon(y))\cap N_\epsilon(y)=\emptyset
\quad\mbox{for all}~\alpha\in\Z^n, j\in\Z.
\end{equation}
\end{lemma}

Using Lemma~\ref{L.fixed}, we show the existence of a set $I\subset\R^n$ satisfying some properties, which will be crucial for the construction of wavelets in $\M_r$.

\begin{proposition}
\label{P.swset}
Let $k\in\Z^n\setminus\{0\}$ and $p\in\Z\setminus\{0\}$. There exists a measurable set $I\subset\R^n$ satisfying the following properties:
\bea
& &|I|>0, \label{E.imp} \\
& & \tau(I)\cap\tau(B^{-p}I)=\emptyset, \label{E.taui} \\
& & d(I)\cap d(I+2B^pk\pi)=\emptyset, \label{E.di} \\
& & I\cup(B^{-p}I+2k\pi) ~\mbox{is contained in an $A$-wavelet set}. \label{E.ikr} 
\eea
\end{proposition}

\proof
Let $Y$ be the collection of fixed points of ${\mathcal F}\cup{\mathcal G}$ as defined in Lemma~\ref{L.fixed}. Notice that the fixed points of $\mathcal F$ is precisely the lattice $2\pi(B^{-p}-I)^{-1}\Z^n$ and those of $\mathcal G$ is the set $\{2\pi(B^{-j}-I)^{-1}k:j\in\Z\setminus\{0\}\}$, which has limit points $\{0,-2k\pi\}$.
Choose any $y\not\in\overline{Y}\cup(2\pi\Z^n)\cup(2\pi B^p\Z^n)$. By Lemma~\ref{L.fixed}, there exists an $\epsilon>0$ such that (\ref{E.faga}) holds. That is, if we take $I=N_\epsilon(y)$, then $I$ satisfies (\ref{E.imp}), (\ref{E.taui}), and (\ref{E.di}). To prove (\ref{E.ikr}) we shall show that $S=I\cup(B^{-p}I+2k\pi)$ satisfies the conditions (a)--(d) of Theorem~\ref{T.subws}. Equations (\ref{E.taui}) and (\ref{E.di}) clearly imply (a) and (b). Since $y\not\in(2\pi\Z^n)\cup(2\pi B^p\Z^n)$, by chosing $\epsilon'<\epsilon$ sufficiently small, if necessary, and taking $I=N_\epsilon'(y)$, we observe that (c) and (d) are also satisfied. Hence, there is an $A$-wavelet set containing 
$I\cup(B^{-p}I+2k\pi)$.
\qed

Now fix an integer $r\geq 0$. Choose $k=k_r\in B^r\Z^n\setminus B^{r+1}\Z^n$ and $p=1$. Then by Proposition~\ref{P.swset}, there is a set $I\subset\R^n$ of positive measure and an $A$-wavelet set $W$ such that
\bea
& & \tau(I)\cap\tau(B^{-1}I)=\emptyset, \label{E.taui1} \\
& & d(I)\cap d(I+2Bk_r\pi)=\emptyset, \label{E.di1} \\
& & I\cup(B^{-1}I+2k_r\pi)\subset W. \label{E.ikr1} 
\eea

Let $J=W\setminus(I\cup(B^{-1}I+2 k_r\pi))$, and define the function $\psi_r$ as follows:
\begin{eqnarray}\label{E.psir}
\hat\psi_r(\xi)=
\begin{cases}
\frac{1}{\sqrt{2}} & 
    {\rm if}~\xi\in I\cup (B^{-1}I)\cup    
    (B^{-1}I+2k_r\pi) \\
-\frac{1}{\sqrt 2} & 
    {\rm if}~\xi\in I+2Bk_r\pi \\
1 & {\rm if}~\xi\in H \\
0 & {\rm otherwise}.
\end{cases}
\end{eqnarray}

\begin{theorem}
\label{T.exp}
For each integer $r\geq 0$, $\psi_r$ is an $A$-wavelet belonging to the equivalence class $\M_r$.
\end{theorem}

\proof
To prove that $\psi_r$ is an $A$-wavelet, we verify that $\psi_r$ satisfies (i)--(iii) of Theorem~\ref{T.wavelet}. Condition (i) is verified once we show that $\sum_{m\in\Z^n}|\hat\psi_r(\xi+2m\pi)|^2=1$ a.e. For, 
\[
\|\hat\psi_r\|^2_2=\int_{\R^n}|\hat\psi_r(\xi)|^2d\xi = \int_{\T^n}\sum_{m\in\Z^n}|\hat\psi_r(\xi+2m\pi)|^2d\xi = (2\pi)^n,
\]
hence, $\|\psi_r\|^2_2=1$. Let $\rho(\xi)=\sum_{m\in\Z^n}|\hat\psi_r(\xi+2k\pi)|^2$. Since $\rho$ is $2\pi\Z^n$-periodic, it suffices to prove that $\rho(\xi)=1$ a.e. on the wavelet set $W$ which is $2\pi\Z^n$-translation equivalent to $\T^n$. Notice that 
\[
W=J\cup I\cup (B^{-1}I+2k_r\pi)~{\rm and}~ 
{\rm supp~}\hat\psi_r=W\cup(B^{-1}I)\cup(I+2Bk_r\pi).
\]
If $\xi\in J$, then $\xi+2m\pi\in {\rm supp~}\hat\psi_r$ if and only if $m=0$, hence $\rho(\xi)=1$. If $\xi\in I$, then $\xi+2m\pi\in {\rm supp~}\hat\psi_r$ if and only if  $m=0$, or $m=Bk_r$. Hence,  $\rho(\xi)=|\hat\psi_r(\xi)|^2+|\hat\psi_r(\xi+2Bk_r\pi)|^2$= 
$(\frac{1}{\sqrt{2}})^2+(-\frac{1}{\sqrt{2}})^2=1$. Finally, if $\xi\in B^{-1}I+2k_r\pi$, then $\xi+2m\pi\in {\rm supp~}\hat\psi_r$ if and only if $m=0$, or $m=-k_r$ so that  $\rho(\xi)=|\hat\psi_r(\xi)|^2+|\hat\psi_r(\xi-2k_r\pi)|^2$= 
$(\frac{1}{\sqrt{2}})^2+(\frac{1}{\sqrt{2}})^2=1$. 

The proof of (ii) is similar.

Now we prove (iii). The term $\hat\psi_r(B^j\xi)\overline{\hat\psi_r(B^j(\xi+2q\pi))}$ is non-zero only when both $B^j\xi$ and $B^j(\xi+2q\pi)$ are in the support of $\hat\psi_r$. Since $B^j(\xi+2q\pi)=B^j\xi+2B^jq\pi$, it is clear from (\ref{E.taui1}), (\ref{E.ikr1}),  and the definition of $\hat\psi_r$ that either $B^jq=k_r$ or $B^jq=Bk_r$. In the first case, we get that $B^j(\xi+2q\pi)\in I+2Bk_r\pi$ and $B^j\xi\in I$ if and only if $B^{j-1}(\xi+2q\pi)\in B^{-1}I+2k_r\pi$ and $B^{j-1}\xi\in B^{-1}I$. Hence, the sum in Theorem~\ref{T.wavelet} (iii) is equal to 
$(\frac{1}{\sqrt{2}})(-\frac{1}{\sqrt{2}})+(\frac{1}{\sqrt{2}})(\frac{1}{\sqrt{2}})=0$. The second case is similar. 

Therefore, by Theorem~\ref{T.wavelet}, it follows that $\psi_r$ is an $A$-wavelet.

Finally, we have to show that $\psi_r\in\M_r$. Recall that 
$E(k)={\rm supp~}\hat\psi_r\cap({\rm supp~}\hat\psi_r+2k\pi)$. From the definition of $\psi_r$, it is clear that $E(k)$ is non-empty if and only if $k=0,\pm k_r,\pm Bk_r$. By our choice, $k_r\in B^r\Z^n\setminus B^{r+1}\Z^n$. Hence, by Theorem~\ref{T.char}(b) and (c), $\psi_r\in\M_r$, $r\geq 0$.
\qed

As we mentioned earlier, it was proved in~\cite{DLS} that there are MSF $A$-wavelets of $L^2(\R^n)$ associated with every dilation matrix. Hence, $\M_\infty$ is non-empty. We just proved that $\M_r$ is non-empty for $r\geq 0$. Therefore, each of the equivalence classes of $A$-wavelets is non-empty.



\begin{thebibliography}{amsplain}
\bibitem{Beh} Behera, B., 
{\it Non-MSF wavelets for the Hardy space $H^2(\R)$}, 
preprint (2002).

\bibitem{BM} Behera, B., and Madan, S., 
{\it Wavelet subspaces invariant under translation operators}, 
preprint (2001).

\bibitem{BS}Bownik, M., and Speegle, D., 
{\it The wavelet dimension function for real dilations and dilations admitting non-MSF wavelets},
Approximation Theory X: Wavelets, Splines, and Applications, 63-85, 
Vanderbilt University Press, 2002. 

\bibitem{Cal} Calogero, A.,
{\it A characterization of wavelets on general lattices}, 
J. Geom. Anal. {\bf 10} (2000), 597--622.

\bibitem{DLS}  Dai, X., Larson, D., and  Speegle, D., 
{\it Wavelet sets in $\R\sp n$}, 
J. Fourier Anal. Appl. {\bf 3} (1997), 451--456.

\bibitem{GM}Gr\"ochenig, K., and  Madych, W.,
{\it Multiresolution analysis, Haar bases, and self-similar tilings of $\R\sp n$},
IEEE Trans. Inform. Theory {\bf 38} (1992), 556--568. 

\bibitem{HW} Hern\'andez, E., and Weiss, G.,
{\sl A First Course on Wavelets}, 
CRC Press, Boca Raton (1996).

\bibitem{Rud} Rudin, W., 
{\sl Fourier analysis on groups}, 
John Wiley and Sons, New York-London (1962).

\bibitem{SW}Schaffer, S., and Weber, E.,
{\it Wavelets with the translation invariance property of order n}, 
preprint (2000).
 
\bibitem{Web} Weber, E., 
{\it On the translation invariance of wavelet subspaces}, 
J. Fourier Anal. Appl. {\bf 6} (2000), 551--558. 

\bibitem{Woj} Wojtaszczyk, P., 
{\sl A mathematical introduction to wavelets}, 
Cambridge University Press, Cambridge (1997).

\end{thebibliography}
\end{document}